# A SIMPLE POLYNOMIAL TIME ALGORITHM TO APPROXIMATE THE PERMANENT WITHIN A SIMPLY EXPONENTIAL FACTOR

ALEXANDER BARVINOK

ABSTRACT. We present a simple randomized polynomial time algorithm to approximate the mixed discriminant of $n$ positive semidefinite $n \times n$ matrices within a factor $2^{O(n)}$. Consequently, the algorithm allows us to approximate in randomized polynomial time the permanent of a given $n \times n$ non-negative matrix within a factor $2^{O(n)}$. When applied to approximating the permanent, the algorithm turns out to be a simple modification of the well-known Godsil-Gutman estimator.

## 1. INTRODUCTION

In this paper we address the question of how to approximate the permanent of a non-negative matrix, and, more generally, the mixed discriminant of positive semidefinite matrices. Our main result is that a simple modification of the well-known Godsil-Gutman estimator ([10], see also Chapter 8 of [19]) yields a randomized polynomial time algorithm, which, given an $n \times n$ non-negative matrix, approximates its permanent within a $2^{O(n)}$ factor. It turns out that the ideas of the algorithm and the proof become more transparent when generalized to mixed discriminants. The first randomized polynomial time algorithm that approximates the permanent within a $2^{O(n)}$ factor was suggested by the author in [5]. The algorithm described in this paper has some advantages compared to the algorithm from [5]. It has a much more transparent structure, it is much easier to implement and it is easily parallelizable. Besides, it sheds some additional light on the properties of the Godsil-Gutman estimator, which was studied in several papers (see [9], [16] and Chapter 8 of [19]).

**(1.1) Permanent.** Let $A = (a_{ij})$ be an $n \times n$ matrix and let $S_n$ be the symmetric group, that is the group of all permutations of the set $\{1, \ldots, n\}$. The number

$$\operatorname{per} A = \sum_{\sigma \in S_n} \prod_{i=1}^{n} a_{i\sigma(i)}$$

is called the *permanent* of $A$. If $A$ is a 0-1 matrix, then per $A$ can be interpreted as the number of perfect matchings in a bipartite graph $G$ on $2n$ vertices $v_1, \ldots, v_n$

---

and $u_1, \ldots, u_n$, where $(v_i, u_j)$ is an edge of $G$ if and only if $a_{ij} = 1$. To compute the permanent of a given 0-1 matrix is a #P-complete problem and even to estimate per $A$ seems to be difficult. Polynomial time algorithms for computing per $A$ are known when $A$ has some special structure, for example, when $A$ is sparse [11] or has small rank [4]. Polynomial time approximation schemes are known for dense 0-1 matrices [12], for "almost all" 0-1 matrices (see [12], [21] and [9]) and for some special 0-1 matrices, such as corresponding to lattice graphs, (see [13] for a survey on approximation algorithms). However, not much is known on how to approximate in polynomial time the permanent of an arbitrary 0-1 matrix (see [15] for the fastest known "mildly exponential" approximation scheme), let alone the permanent of an arbitrary non-negative matrix.

Let $t_1, \ldots, t_n$ be real variables. Then the permanent per $A$ can be expressed as the coefficient of $t_1 \cdots t_n$ in the product of linear forms:

$$(1.1.1) \qquad \text{per } A = \frac{\partial^n}{\partial t_1 \ldots \partial t_n} \prod_{i=1}^{n} \sum_{j=1}^{n} a_{ij} t_j.$$

**(1.2) Mixed Discriminant.** Let $Q_1, \ldots, Q_n$ be $n \times n$ symmetric matrices and let $t_1, \ldots, t_n$ be real variables. Then $\det(t_1 Q_1 + \ldots + t_n Q_n)$ is a homogeneous polynomial of degree $n$ in $t_1, \ldots, t_n$. The number

$$D(Q_1, \ldots, Q_n) = \frac{\partial^n}{\partial t_1 \ldots \partial t_n} \det(t_1 Q_1 + \ldots + t_n Q_n)$$

is called the mixed discriminant of $Q_1, \ldots, Q_n$. Sometimes the normalizing factor $1/n!$ is used (cf. [18]). The mixed discriminant $D(Q_1, \ldots, Q_n)$ is a polynomial in the entries of $Q_1, \ldots, Q_n$ with coefficients $-1$ and $1$.

The mixed discriminant can be considered as a generalization of the permanent. Indeed, from (1.1.1) we deduce that for diagonal matrices $Q_1, \ldots, Q_n$

$(1.2.1)$
$$D(Q_1, \ldots, Q_n) = \text{per } A, \quad \text{where} \quad Q_i = \text{diag}\{a_{i1}, \ldots, a_{in}\} \quad \text{and} \quad A = (a_{ij}).$$

Mixed discriminants were introduced by A.D. Aleksandrov in his proof of the Aleksandrov - Fenchel inequality for mixed volumes ([2], see also [18]). The relation between the mixed discriminant and the permanent was used in the proof of the Van der Waerden conjecture for permanents of doubly stochastic matrices (see [6]).

The mixed discriminant is linear in each argument, that is

$$D(Q_1, \ldots, Q_{i-1}, \alpha Q_i + \beta Q_i', Q_{i+1}, \ldots, Q_n)$$

$$= \alpha D(Q_1, \ldots, Q_{i-1}, Q_i, Q_{i+1}, \ldots, Q_n) + \beta D(Q_1, \ldots, Q_{i-1}, Q_i', Q_{i+1}, \ldots, Q_n),$$

and $D(Q_1, \ldots, Q_n) \geq 0$ provided $Q_1, \ldots, Q_n$ are positive semidefinite (see [18]).

The paper is organized as follows. In Section 2, we describe the algorithms for computing mixed discriminants and permanents, state the main results, and discuss them. Section 3 addresses the behavior of quadratic forms on random vectors drawn from the Gaussian distribution in $\mathbb{R}^n$ – our main tool for analyzing the algorithms. Section 4 contains proofs of the main results of the paper. In Section 5, we present an application of the mixed discriminant for counting problems and also make some remarks on "binary versions" of our algorithms.

## 2. Main Results

Our algorithms use two standard procedures from linear algebra: computing the determinant of a matrix and computing a decomposition $Q = TT^*$, where $Q$ is a positive semidefinite matrix and $T^*$ is the transpose of $T$. One can compute the determinant of a given $n \times n$ matrix using $O(n^3)$ arithmetic operations (see, for example, Chapter 2, Section 7 of [8]). For a positive semidefinite $n \times n$ matrix $Q$ one can compute a decomposition $Q = TT^*$, where $T$ is a lower triangular matrix, by using $O(n^3)$ arithmetic operations and $n$ times taking square root from a non-negative number (see, for example, Chapter 2, Section 10 of [8]). The "random" part of the algorithms consists of sampling $n^2$ random variables independently from the standard Gaussian distribution in $\mathbb{R}$ with the density

$$\psi_1(x) = \frac{1}{\sqrt{2\pi}} e^{-x^2/2}.$$

Various ways to simulate this distribution from the uniform distribution on the interval $[0,1]$ are described in Section 3.4.1 (C) of [17]. In particular, this allows us to sample vectors $x = (x_1, \ldots, x_n) \in \mathbb{R}^n$ from the standard $n$-dimensional Gaussian distribution in $\mathbb{R}^n$ with the density

$$\psi_n(x) = (2\pi)^{-n/2} \exp\{-\|x\|^2/2\}, \quad \text{where} \quad \|x\|^2 = x_1^2 + \ldots + x_n^2,$$

by sampling the coordinates $x_1, \ldots, x_n$ independently from the one-dimensional Gaussian distribution. We will write $\psi(x)$ instead of $\psi_n(x)$ if the choice of the ambient space $\mathbb{R}^n$ is clear from the context.

To simplify our analysis, we assume that we operate with real numbers and that we can perform arithmetic operations, take square root and sample from the Gaussian distribution exactly. In Section 5.1 we briefly discuss how to adjust our algorithms to the "binary model" of computation, that is, if we allow only integers (rational numbers) and arithmetic operations on them. As a general source on algorithms and complexity we use [1].

Now we describe our main algorithm.

**(2.1) Algorithm for computing the mixed discriminant.**

**Input:** Positive semidefinite $n \times n$ matrices $Q_1, \ldots, Q_n$.

**Output:** A number $\alpha$ approximating the mixed discriminant $D(Q_1, \ldots, Q_n)$.

**Algorithm:** For $i = 1, \ldots, n$ compute a decomposition $Q_i = T_i T_i^*$. Sample independently $n$ vectors $u_1, \ldots, u_n$ at random from the standard Gaussian distribution in $\mathbb{R}^n$ with the density $\psi(x) = (2\pi)^{-n/2} \exp\{-\|x\|^2/2\}$. Compute

$$\alpha = \Big( \det[T_1 u_1, \ldots, T_n u_n] \Big)^2,$$

the squared determinant of the matrix with the columns $T_1 u_1, \ldots, T_n u_n$.
Output $\alpha$.

**(2.2) Theorem.**
(2.2.1) *The expectation of $\alpha$ is the mixed discriminant $D(Q_1, \ldots, Q_n)$;*
(2.2.2) *For any $C > 1$ the probability that*

$$\alpha \geq C \cdot D(Q_1, \ldots, Q_n)$$

*does not exceed $C^{-1}$;*
(2.2.3) *Let*

$$c_0 = \exp\left\{ \frac{4}{\sqrt{2\pi}} \int_0^{+\infty} (\ln t) e^{-t^2/2} \, dt \right\} \approx 0.2807297419.$$

*Then for any $1 > \epsilon > 0$ the probability that*

$$\alpha \leq (\epsilon c_0)^n D(Q_1, \ldots, Q_n)$$

*does not exceed $\dfrac{8}{n \ln^2 \epsilon}$.*

The algorithm becomes especially simple when we apply it to compute the mixed discriminant of diagonal matrices, that is the permanent of a non-negative matrix.

**(2.3) Algorithm for computing the permanent of a non-negative matrix.**
**Input:** Non-negative $n \times n$ matrix $A$.
**Output:** A number $\alpha$ approximating the permanent per $A$.
**Algorithm:** Sample independently $n^2$ numbers $u_{ij} : i, j = 1, \ldots, n$ at random from the standard Gaussian distribution in $\mathbb{R}$ with the density $\psi(x) = \dfrac{1}{\sqrt{2\pi}} e^{-x^2/2}$. Compute the $n \times n$ matrix $B = (b_{ij})$, where $b_{ij} = u_{ij} \sqrt{a_{ij}}$. Compute $\alpha = (\det B)^2$. Output $\alpha$.

**(2.4) Theorem.**
(2.4.1) *The expectation of $\alpha$ is the permanent per $A$;*
(2.4.2) *For any $C > 1$ the probability that*

$$\alpha \geq C \cdot \text{per } A$$

*does not exceed $C^{-1}$;*

(2.4.3) *Let*

$$c_0 = \exp\left\{\frac{4}{\sqrt{2\pi}} \int_0^{+\infty} (\ln t) e^{-t^2/2} \, dt\right\} \approx 0.2807297419.$$

*Then for any $1 > \epsilon > 0$ the probability that*

$$\alpha \leq (\epsilon c_0)^n \mathrm{per}\ A$$

*does not exceed* $\dfrac{8}{n \ln^2 \epsilon}$.

As implied by (1.2.1), Theorem 2.4 is a straightforward corollary of Theorem 2.2. From (2.2.2) and (2.2.3) it follows that for all sufficiently large $n$, the output $\alpha$ of Algorithm 2.1 satisfies the inequalities

$$(0.28)^n D(Q_1, \ldots, Q_n) \leq \alpha \leq 3D(Q_1, \ldots, Q_n)$$

with probability at least 0.6. We can make the probability as high as $1 - \epsilon$ for any $\epsilon > 0$ by running the algorithm independently $O\left(\ln(\epsilon^{-1})\right)$ times and taking the median of the computed $\alpha$'s (cf. [14]). Hence we get a randomized polynomial time algorithm approximating the mixed discriminant of positive semidefinite matrices (and hence the permanent of a non-negative matrix) within a simply exponential factor $2^{O(n)}$.

**(2.5) Relation to the Godsil-Gutman estimator.** It is seen immediately that Algorithm 2.3 is a modification of the Godsil-Gutman estimator (see [10]) and Chapter 8 of [19]). Indeed, in the Godsil-Gutman estimator we sample $u_{ij}$ from the binary distribution:

$$u_{ij} = \begin{cases} 1 & \text{with probability } 1/2, \\ -1 & \text{with probability } 1/2. \end{cases}$$

Furthermore, parts (2.4.1) and (2.4.2) of Theorem 2.4 remain true as long as we sample $u_{ij}$ independently from some distribution with the expectation 0 and variance 1. However, (2.4.3) is not true for the binary distribution. Indeed, let $A = (a_{ij})$ be the following $n \times n$ matrix

$$a_{ij} = \begin{cases} 1 & \text{if } i = j \text{ or } \{i, j\} = \{1, 2\}, \\ 0 & \text{elsewhere}. \end{cases}$$

Then per $A = 2$. If $u_{ij} \in \{-1, 1\}$ are sampled from the binary distribution, then $\alpha = (u_{11}u_{22} - u_{12}u_{21})^2$. So, we get that $\alpha = 0$ with probability $1/2$ and $\alpha = 4$ with probability $1/2$. Apparently, sampling from continuous distributions allows us to approximate better.

Another (discrete) distribution for the Godsil-Gutman estimator was studied in [16] and [9]. It was shown that for a "typical" 0-1 matrix we get a polynomial approximation scheme for the permanent [9], whereas for a "worst case" 0-1 matrix it allows us to construct an approximation scheme [16] whose complexity, while still exponential, is significantly better than the complexity of known exact methods (cf. Chapter 7 of [20]).

**(2.6) How well can we approximate the permanent in polynomial time?**
Suppose we have a polynomial time (probabilistic or deterministic) algorithm that for any given $n \times n$ (non-negative or 0-1) matrix $A$ computes a number $\alpha$ such that $\phi(n)\text{per } A \leq \alpha \leq \text{per } A$. What sort of function might $\phi(n)$ be? For an $n \times n$ matrix $A$ and $k > 0$ let us construct the $nk \times nk$ block-diagonal matrix $A \otimes I_k$, having $k$ diagonal copies of $A$. We observe that $A \otimes I_k$ is non-negative if $A$ is non-negative and $A \otimes I_k$ is 0-1 if $A$ is 0-1, and that per $A = \left(\text{per } A \otimes I_k\right)^{1/k}$. Applying our algorithm to $A \otimes I_k$ and taking the root we get an approximation of per $A$ within a factor $\phi^{1/k}(nk)$. So, if we suppose that $\phi$ is the best possible, we may assume that $\phi(n) \geq \phi^{1/k}(nk)$. There are few reasonable choices for such functions $\phi$.

a) $\phi(n) \equiv 1$. This doesn't look likely, given that the problem is #P-hard.

b) For $\epsilon > 0$ one can choose $\phi_\epsilon(n) = 1 - \epsilon$ and the algorithm is polynomial in $\epsilon^{-1}$. In the author's opinion this conjecture is overly optimistic, especially for arbitrary non-negative matrices.

c) For $\epsilon > 0$ one can choose $\phi_\epsilon(n) = (1-\epsilon)^n$, but the algorithm is *not* polynomial in $\epsilon^{-1}$. This type of approximation was conjectured by V.D. Milman.

d) $\phi(n) = c^n$ for some fixed constant $c$. This is the type of a bound achieved by Algorithm 2.3.

We note that functions $\phi(n)$, decreasing faster than a simply exponential function of $n$ (for example, $\phi(n) = 1/n!$) are not interesting since they are beaten by $c^n$ of Algorithm 2.3. The author does not know if the constant $c \approx 0.28$ from Theorem 2.4 can be improved.

## 3. The Gaussian Measure and Quadratic Forms

**(3.1) Notation.**    Let us fix the standard orthonormal basis in $\mathbb{R}^n$. Thus we can identify $n \times n$ matrices with linear operators on $\mathbb{R}^n$.

Let $\langle \cdot, \cdot \rangle$ be the standard inner product on $\mathbb{R}^n$, that is

$$\langle x, y \rangle = x_1 y_1 + \ldots + x_n y_n, \quad \text{where } x = (x_1, \ldots, x_n) \text{ and } y = (y_1, \ldots, y_n).$$

We denote the corresponding Euclidean norm by $\| \cdot \|$.

Let $\mu$ be the standard Gaussian measure on $\mathbb{R}^n$ with the density

$$\psi(x) = (2\pi)^{-n/2} \exp\{-\|x\|^2/2\}.$$

For a $\mu$-integrable function $f : \mathbb{R}^n \longrightarrow \mathbb{R}$ we define its expectation

$$\mathbf{E}(f) = \int_{\mathbb{R}^n} f \; d\mu = \int_{\mathbb{R}^n} f(x)\psi(x) \; dx.$$

Furthermore, if $F : \mathbb{R}^n \longrightarrow \mathbb{R}^m$, $F(x) = (f_1(x), \ldots, f_m(x))$, we let $\mathbf{E}(F) = \left(\mathbf{E}(f_1), \ldots, \mathbf{E}(f_m)\right) \in \mathbb{R}^m$.

If $x \in \mathbb{R}^n$ is a vector then we denote by $x \otimes x$ the $n \times n$ matrix whose $(i,j)$-th entry is $x_i \cdot x_j$. So, $x \otimes x$ is a positive semidefinite matrix of rank 1, provided $x \neq 0$.

Let $Q$ be a symmetric $n \times n$ matrix and $q(x) = \langle x, Qx \rangle$ be the corresponding quadratic form. We define the trace of $q$:

$$\operatorname{tr}(q) = \sum_{i=1}^{n} Q_{ii}.$$

We start with a simple lemma.

**(3.2) Lemma.**
(3.2.1) *Let $q : \mathbb{R}^n \longrightarrow \mathbb{R}$ be a quadratic form. Then*

$$\mathbf{E}(q) = \operatorname{tr}(q).$$

(3.2.2) *Let us fix an $n \times n$ matrix $T$. Then*

$$\mathbf{E}(Tu \otimes Tu) = TT^*,$$

*where $u \in \mathbb{R}^n$ is sampled from the standard Gaussian distribution $\mu$.*

*Proof.* Since $\mathbf{E}(x_i^2) = 1$ for $i = 1, \ldots n$ and $\mathbf{E}(x_i x_j) = 0$, the proof of (3.2.1) follows. Therefore, $\mathbf{E}(u \otimes u) = I$, the identity matrix, and hence

$$\mathbf{E}(Tu \otimes Tu) = \mathbf{E}\big(T(u \otimes u)T^*\big) = T\mathbf{E}(u \otimes u)T^* = TT^*,$$

and (3.2.2) follows. $\qquad\qquad\qquad\qquad\qquad\qquad\qquad\qquad\qquad\qquad\qquad\qquad\square$

Now we prove the main result of this section.

**(3.3) Theorem.** *Let $q : \mathbb{R}^n \longrightarrow \mathbb{R}$ be a positive semidefinite quadratic form, such that $\mathbf{E}(q) = 1$. Let*

$$C_0 = \frac{4}{\sqrt{2\pi}} \int_0^{+\infty} (\ln t) e^{-t^2/2} \, dt \approx -1.270362845.$$

*Then*

(3.3.1) $$C_0 \leq \mathbf{E}(\ln q) \leq 0$$

*and*

(3.3.2) $$0 \leq \mathbf{E}(\ln^2 q) \leq 8.$$

*Proof.* Since $\ln$ is a concave function, we have $\mathbf{E}(\ln q) \leq \ln(\mathbf{E}(q)) = 0$. Let us decompose $q$ into a non-negative linear combination $q = \lambda_1 q_1 + \ldots + \lambda_n q_n$ of positive semidefinite forms $q_i$ of rank 1. We can scale $q_i$ so that $\mathbf{E}(q_i) = 1$ for $i = 1, \ldots, n$ and then we have $\lambda_1 + \ldots + \lambda_n = 1$. In fact, one can choose $\lambda_i$ to be the eigenvalues of $q$ and $q_i = \langle x, u_i \rangle^2$, where $u_i$ is the corresponding unit eigenvector. Since $\ln$ is a concave function we have $\ln(\lambda_1 q_1 + \ldots + \lambda_n q_n) \geq \lambda_1 \ln q_1 + \ldots + \lambda_n \ln q_n$. Furthermore, if $q_i$ is a positive semidefinite form of rank 1 such that $\mathbf{E}(q_i) = 1$, then by an orthogonal transformation of the coordinates it can be brought into the form $q_i(x) = x_1^2$ (see (3.2.1)). Therefore, $\mathbf{E}(\ln q_i) = \mathbf{E}(\ln x_1^2) = C_0$ and

$$\mathbf{E}(\ln q) \geq \lambda_1 \mathbf{E}(\ln q_1) + \ldots + \lambda_n \mathbf{E}(\ln q_n) = (\lambda_1 + \ldots + \lambda_n)C_0 = C_0,$$

so (3.3.1) is proven (we note that this reasoning proves that $\mathbf{E}(\ln q)$ is well-defined).

Let $X = \left\{ x \in \mathbb{R}^n : q(x) \leq 1 \right\}$ and $Y = \mathbb{R}^n \setminus X$. Then

$$\mathbf{E}(\ln^2 q) = \int_X \psi(x) \ln^2 q(x) \; dx + \int_Y \psi(x) \ln^2 q(x) \; dx.$$

Let us estimate the first integral. Decomposing $q = \lambda_1 q_1 + \ldots + \lambda_n q_n$ as above, we get $\ln q \geq \lambda_1 \ln q_1 + \ldots + \lambda_n \ln q_n$. Since $\ln q(x) \leq 0$ for $x \in X$, we get that

$$\ln^2 q(x) \leq \sum_{i,j=1}^n \lambda_i \lambda_j \big(\ln q_i(x)\big)\big(\ln q_j(x)\big)$$

for $x \in X$. Therefore,

$$\int_X \psi(x) \ln^2 q(x) \; dx \leq \sum_{i,j=1}^n \lambda_i \lambda_j \int_X \psi(x) \big(\ln q_i(x)\big)\big(\ln q_j(x)\big) \; dx$$

$$\leq \sum_{i,j=1}^n \lambda_i \lambda_j \left( \int_X \psi(x) \ln^2 q_i(x) \; dx \right)^{1/2} \left( \int_X \psi(x) \ln^2 q_j(x) \; dx \right)^{1/2}$$

(we applied the Cauchy-Schwartz inequality)

$$\leq \sum_{i,j=1}^n \lambda_i \lambda_j \Big(\mathbf{E}(\ln^2 q_i)\Big)^{1/2} \Big(\mathbf{E}(\ln^2 q_j)\Big)^{1/2}.$$

Now, as in the proof of (3.3.1) we have

$$\mathbf{E}(\ln^2 q_i) = \mathbf{E}(\ln^2 x_1^2) = \frac{8}{\sqrt{2\pi}} \int_0^{+\infty} (\ln^2 t) e^{-t^2/2} \; dt \approx 6.548623960 \leq 7.$$

Summarizing, we get

$$\int_X \psi(x) \ln^2 q(x) \; dx \leq \sum_{i,j=1}^n \lambda_i \lambda_j \Big(\mathbf{E}(\ln^2 q_i)\Big)^{1/2} \Big(\mathbf{E}(\ln^2 q_j)\Big)^{1/2} \leq 7 \sum_{i,j=1}^n \lambda_i \lambda_j = 7.$$

Since for $0 \leq \ln t \leq \sqrt{t}$ for $t \geq 1$ we have

$$\int_Y \psi(x) \ln^2 q(x) \; dx \leq \int_Y q(x) \psi(x) \; dx \leq \mathbf{E}(q) = 1.$$

Therefore, $\mathbf{E}(\ln^2 q) \leq 7 + 1 = 8$ and (3.3.2) is proven. $\qquad\square$

Finally, we will need the following simple result.

**(3.4) Lemma.** *Let $u_1, \ldots, u_n$ be vectors from $\mathbb{R}^n$. Then*

$$D(u_1 \otimes u_1, \ldots, u_n \otimes u_n) = \Big(\det[u_1, \ldots, u_n]\Big)^2,$$

*the squared determinant of the matrix with the columns $u_1, \ldots, u_n$.*

*Proof.* Let $e_1, \ldots, e_n$ be the standard orthonormal basis of $\mathbb{R}^n$. Let $G$ be the operator such that $G(e_i) = u_i$ for $i = 1, \ldots, n$. Then $u_i \otimes u_i = G(e_i \otimes e_i)G^*$ and from the definition of Section 1.2 we get

$$D(u_1 \otimes u_1, \ldots, u_n \otimes u_n) = \frac{\partial^n}{\partial t_1 \ldots \partial t_n} \det\Big(t_1 u_1 \otimes u_1 + \ldots + t_n u_n \otimes e_n\Big)$$

$$= \frac{\partial^n}{\partial t_1 \ldots \partial t_n} \det\Big(G(t_1 e_1 \otimes e_1 + \ldots + t_n e_n \otimes e_n)G^*\Big)$$

$$= \det(GG^*)\frac{\partial^n}{\partial t_1 \ldots \partial t_n} \det\Big(t_1 e_1 \otimes e_1 + \ldots + t_n e_n \otimes e_n\Big) = (\det G)^2.$$

$\square$

## 4. PROOF OF THE MAIN RESULTS

**(4.1) Conditional expectations.** In this subsection, we summarize some general results on measures and integration, which we exploit later. As a general source we use [3].

Suppose that we have $k$ copies of the Euclidean space $\mathbb{R}^n$, each with the standard Gaussian probability measure $\mu$. We will consider functions $f: \mathbb{R}^n \times \ldots \times \mathbb{R}^n \longrightarrow \mathbb{R}$ that are defined almost everywhere and integrable with respect to the measure $\nu_k = \mu \times \ldots \times \mu$. Let $f(u_1, \ldots, u_k)$ be such a function. Then for almost all $(k-1)$-tuples $(u_1, \ldots, u_{k-1})$ with $u_i \in \mathbb{R}^n$, the function $f(u_1, \ldots, u_{k-1}, \cdot)$ is integrable (Fubini's Theorem) and we can define the *conditional expectation* $g(u_1, \ldots, u_{k-1}) = \mathbf{E}_k(f)$ by letting

$$\mathbf{E}_k(f)(u_1, \ldots, u_{k-1}) = \int_{\mathbb{R}^n} f(u_1, \ldots, u_{k-1}, u_k)\psi(u_k) \, du_k.$$

Fubini's theorem implies that

$$\mathbf{E}(f) = \mathbf{E}_1 \ldots \mathbf{E}_k(f),$$

where $\mathbf{E}$ is the expectation with respect to the product measure $\nu_k$. Tonelli's Theorem states that if $f$ is $\nu_k$-measurable and non-negative almost surely with respect to $\nu_k$ then $f$ is $\nu_k$-integrable, provided $\mathbf{E}_1 \ldots \mathbf{E}_k(f) < +\infty$.

If $f(u_1, \ldots, u_i)$ is a function of $i < k$ arguments, we may formally extend it to $\mathbb{R}^n \times \ldots \times \mathbb{R}^n$ ($k$ times) by letting $f(u_1, \ldots, u_k) = f(u_1, \ldots, u_i)$. The distribution of values of $f(u_1, \ldots, u_i)$ with respect to $\nu_i$ is the same as the distribution of values of $f(u_1, \ldots, u_k)$ with respect to $\nu_k$. In particular, if $f(u_1, \ldots, u_i)$ is $\nu_i$-integrable then $f(u_1, \ldots, u_k)$ is $\nu_k$-integrable.

We note the following useful facts:

(4.1.1) The linear operator $\mathbf{E}_k$ is *monotone*, that is, if $f(u_1, \ldots, u_k) \leq g(u_1, \ldots, u_k)$ almost surely with respect to $\nu_k$, then $\mathbf{E}_k(f) \leq \mathbf{E}_k(g)$ almost surely with respect to $\nu_{k-1}$;

(4.1.2) If $f(u_1, \ldots, u_k)$ is integrable and $g(u_1, \ldots, u_i)$, $i < k$ is a function, then $\mathbf{E}_k(gf) = g\mathbf{E}_k(f)$;

(4.1.3) If $f = a$ is a constant almost surely with respect to $\nu_k$, then $\mathbf{E}_k(f) = a$ almost surely with respect to $\nu_{k-1}$.

First, we prove the following technical lemma (a martingale inequality).

**(4.2) Lemma.** *Let $f_k(u_1, \ldots, u_k)$, $k = 1, \ldots, n$ be integrable functions on $\mathbb{R}^n \times \ldots \times \mathbb{R}^n$ ($k$ times), and let $\nu_n = \mu \times \ldots \times \mu$ ($n$ times). Suppose that for some numbers $a$ and $b$ we have*

$$a \leq \mathbf{E}_k(f_k) \quad and \quad \mathbf{E}_k(f_k^2) \leq b \qquad almost\ surely\ with\ respect\ to\ \nu_{k-1}$$

*for $k = 1, \ldots, n$.*

*Then for any $\delta > 0$ we have*

$$\nu_n\Big\{(u_1, \ldots, u_n) : \ \frac{1}{n}\sum_{k=1}^n f_k(u_1, \ldots, u_k) \leq a - \delta\Big\} \leq \frac{b}{n\delta^2}.$$

*Proof.* Let $g_k = \mathbf{E}_k(f_k)$ and $h_k = f_k - g_k$. Since $g_k$ does not depend on $u_k$, using (4.1.2) we have $\mathbf{E}_k(h_k^2) = \mathbf{E}_k(f_k^2) - 2\mathbf{E}_k(g_k f_k) + \mathbf{E}_k(g_k^2) = \mathbf{E}_k(f_k^2) - g_k^2$. Summarizing, we may write

$$f_k = g_k + h_k, \quad \text{where} \quad \mathbf{E}_k(h_k) = 0, \quad g_k \geq a, \qquad \text{and} \qquad \mathbf{E}_k(h_k^2) \leq b$$

almost surely with respect to $\nu_{k-1}$. Let

$$H(u_1, \ldots, u_n) = \frac{1}{n}\sum_{k=1}^n h_k(u_1, \ldots, u_k).$$

Now for $U = (u_1, \ldots, u_n)$ we have

$$\nu_n\Big\{U : \frac{1}{n}\sum_{k=1}^n f_k(u_1, \ldots, u_k) \leq a - \delta\Big\}$$

$$= \nu_n\Big\{U : \ H(U) + \frac{1}{n}\sum_{k=1}^n g_k(u_1, \ldots, u_{k-1}) \leq a - \delta\Big\}$$

$$\leq \nu_n\Big\{U : \ H(U) \leq -\delta\Big\} \leq \frac{\mathbf{E}(H^2)}{\delta^2} = \frac{1}{\delta n^2}\sum_{k=1}^n \mathbf{E}(h_k^2) + \frac{2}{\delta n^2}\sum_{1 \leq i < j \leq n} \mathbf{E}(h_i h_j)$$

(we used Chebyshev's inequality).

We note that it is legitimate to pass to global expectations $\mathbf{E}$ here. Indeed, since $h_k^2$ is non-negative and $\mathbf{E}_k h_k^2 \leq \mathbf{E}_k f_k^2 \leq b$ it follows by (4.1.1), (4.1.3) and Tonelli's Theorem that $h_k^2$ is $\nu_n$-integrable. Since $|h_i h_j| \leq (h_i^2 + h_j^2)/2$, the products $h_i h_j$ are also $\nu_n$-integrable. Therefore, $H^2$ is $\nu_n$-integrable.

Since $h_k$ does not depend on $u_{k+1}, \ldots, u_n$, using (4.1.2) we have $\mathbf{E}(h_k^2) = \mathbf{E}_1 \ldots \mathbf{E}_n h_k^2 = \mathbf{E}_1 \ldots \mathbf{E}_k h_k^2$ and since $\mathbf{E}_k h_k^2 \leq b$ almost surely on $\nu_{k-1}$, by (4.1.1) and (4.1.3) we get that $\mathbf{E}(h_k^2) \leq b$ for each $k = 1, \ldots, n$. Furthermore, since $h_i$ does not depend on $u_{i+1}, \ldots, u_n$ for $j > i$, using (4.1.2) and (4.1.3) we have

$$\mathbf{E}(h_i h_j) = \mathbf{E}_1 \ldots \mathbf{E}_n(h_i h_j) = \mathbf{E}_1 \ldots \mathbf{E}_j(h_i h_j) = \mathbf{E}_i \ldots \mathbf{E}_{j-1} h_i \mathbf{E}_j(h_j) = 0.$$

The proof now follows. □

Now we are ready to prove the main results of the paper.

**(4.3) Proof of Theorem 2.2.** The output $\alpha$ of Algorithm 2.1 is a function of vectors $u_1, \ldots, u_n \in \mathbb{R}^n$, which are drawn at random from the standard Gaussian distribution $\mu$ in $\mathbb{R}^n$. We consider the distribution of $\alpha(u_1, \ldots, u_n)$ with respect to the product measure $\nu_n = \mu \times \ldots \times \mu$ ($n$ times).

Using Lemma 3.4, we may write

$$\alpha = \alpha(u_1, \ldots, u_n) = \Big(\det[T_1 u_1, \ldots, T_n u_n]\Big)^2 = D(T_1 u_1 \otimes T_1 u_1, \ldots, T_n u_n \otimes T_n u_n).$$

For $k = 0, \ldots, n$ let

$$\alpha_k(u_1, \ldots, u_k) = D(T_1 u_1 \otimes T_1 u_1, \ldots, T_k u_k \otimes T_k u_k, Q_{k+1}, \ldots, Q_n).$$

In particular, $\alpha_0 = D(Q_1, \ldots, Q_n)$ and $\alpha_n = \alpha$. Since the mixed discriminant of positive semidefinite matrices is non-negative (Section 1.2), we deduce that $\alpha_k(u_1, \ldots, u_k)$ is non-negative. By (3.2.2) we have $\mathbf{E}_k(T_k u_k \otimes T_k u_k) = T_k T_k^* = Q_k$. Since $\alpha_k$ is a polynomial in the coordinates of $u_1, \ldots, u_k$ and the mixed discriminant is linear in every argument (see Section 1.2) we can interchange $D$ and $\mathbf{E}_k$:

$$\mathbf{E}_k(\alpha_k)$$
$$= D\big(T_1 u_1 \otimes T_1 u_1, \ldots, T_{k-1} u_{k-1} \otimes T_{k-1} u_{k-1}, \mathbf{E}_k(T_k u_k \otimes T_k u_k), Q_{k+1}, \ldots, Q_n\big)$$
$$= D\big(T_1 u_1 \otimes T_1 u_1, \ldots, T_{k-1} u_{k-1} \otimes T_{k-1} u_{k-1}, Q_k, Q_{k+1}, \ldots, Q_n\big) = \alpha_{k-1}.$$

Since $\alpha_k$ is non-negative and $\nu_n$-measurable, applying theorems of Tonelli and Fubini, we have

$$(4.3.1) \qquad \mathbf{E}(\alpha) = \mathbf{E}_1 \mathbf{E}_2 \ldots \mathbf{E}_n(\alpha_n) = \mathbf{E}_1 \ldots \mathbf{E}_k(\alpha_k) = \alpha_0 = D(Q_1, \ldots, Q_n).$$

and (2.2.1) is proven. Since $\alpha(u_1, \ldots, u_n)$ is non-negative, (2.2.2) follows by Chebyshev's inequality.

If $D(Q_1, \ldots, Q_n) = 0$ then by (2.2.1) and non-negativity of $\alpha$ it follows that $\alpha$ is identically 0 and (2.2.3) would follow. Therefore, without loss of generality, we may assume that $D(Q_1, \ldots, Q_n) > 0$. Since $\alpha_k(u_1, \ldots, u_k)$ is a non-negative polynomial, by (4.3.1) we conclude that $\alpha_k(u_1, \dots, u_k) > 0$ almost surely with respect to $\nu_k$. For $k = 1, \dots, n$ let

$$g_k(u_1, \dots, u_k) = \frac{\alpha_k(u_1, \dots, u_k)}{\alpha_{k-1}(u_1, \dots, u_{k-1})} = \frac{\alpha_k}{\mathbf{E}_k(\alpha_k)},$$

Hence we may write

$$\frac{\alpha(u_1, \dots, u_n)}{D(Q_1, \dots, Q_n)} = \prod_{k=1}^{n} g_k(u_1, \dots, u_k)$$

almost surely with respect to $\nu_n$. Since the mixed discriminant is linear in each argument, $g_k(u_1, \dots, u_k)$ is a quadratic form in $u_k$ for any fixed $u_1, \dots, u_{k-1}$, such that $\alpha_{k-1}(u_1, \dots, u_{k-1}) > 0$. Furthermore, since the mixed discriminant is non-negative for positive semidefinite arguments, we conclude that $g_k(u_1, \dots, u_k)$ is a positive semidefinite quadratic form in $u_k$ for every such choice of $u_1, \dots, u_{k-1}$. Since $\mathbf{E}_k(g_k) = 1$ almost surely with respect to $\nu_{k-1}$, by Theorem 3.3 we conclude that

$$C_0 \le \mathbf{E}_k(\ln g_k) \le 0 \quad \text{with} \quad C_0 = \frac{4}{\sqrt{2\pi}} \int_0^{+\infty} (\ln t) e^{-t^2/2} \, dt$$

and that $|\mathbf{E}_k(\ln^2 g_k)| \le 8$ almost surely with respect to $\nu_{k-1}$. In particular, since $\ln^2 g_k$ is non-negative almost surely with respect to $\nu_k$, we deduce that $\ln^2 g_k$ is $\nu_n$-integrable, and since $|\ln g_k| \le 1 + \ln^2 g_k$ we deduce that $\ln g_k$ is $\nu_n$-integrable.

Now for $U = (u_1, \dots, u_n)$ we have

$$\nu_n \left\{ U : \frac{\alpha(u_1, \dots, u_n)}{D(Q_1, \dots, Q_n)} \le (\epsilon c_0)^n \right\} = \nu_n \left\{ U : \frac{1}{n} \ln \left( \frac{\alpha(u_1, \dots, u_n)}{D(Q_1, \dots, Q_n)} \right) \le C_0 + \ln \epsilon \right\}$$

$$= \nu_n \left\{ (u_1, \dots, u_n) : \frac{1}{n} \sum_{k=1}^{n} \ln g_k(u_1, \dots, u_k) \le C_0 + \ln \epsilon \right\}.$$

To complete the proof of (2.2.3), we use Lemma 4.2 with $f_k = \ln g_k$, $a = C_0$, $b = 8$ and $\delta = -\ln \epsilon$.                                                               $\square$

**(4.4) Proof of Theorem 2.4.** For $i = 1, \dots, n$ let us define the diagonal matrices $Q_i = \mathrm{diag}\{a_{i1}, \dots, a_{in}\}$. Algorithm 2.3 with the input $A$ is Algorithm 2.1 with the input $Q_1, \dots, Q_n$. By (1.2.1) per $A = D(Q_1, \dots, Q_n)$, and the proof follows by Theorem 2.2.                                                               $\square$

## 5. REMARKS

**(5.1) The algorithms in the binary model.** Suppose we want to operate in the standard binary (Turing) model of computation. That is, we allow arithmetic operations with integral (rational) numbers stored as bit strings (see [1]). In the probabilistic setting, we suppose also that we can generate a random bit, that is we can sample a random variable $x \in \{0, 1\}$, which assumes either value with probability $1/2$. Algorithms 2.1 and 2.3 can be transformed into randomized polynomial time algorithms in this model, which approximate the mixed discriminant, resp. the permanent within a $2^{O(n)}$ factor. This transformation is relatively straightforward for Algorithm 2.3 and first we briefly sketch it here, omitting tedious details (although it may not be the most efficient binary version).

We are given a non-negative integer matrix $A$ and we want to approximate per $A$. To compute the matrix $B$ from Algorithm 2.3 we need to compute square roots $\sqrt{a_{ij}}$ and to sample variables $u_{ij}$ from the standard Gaussian distribution in $\mathbb{R}$. It is known that the square root of a positive rational number can be computed within any given error $\epsilon > 0$ in time that is polynomial in $\ln \epsilon^{-1}$. We observe that if we choose $u_{ij} : 1 \leq i, j \leq n$ independently from the standard Gaussian distribution in $\mathbb{R}$, then we will have $|u_{ij}| \leq n$ for all $i, j$ with the probability that goes to 1 exponentially fast as $n \longrightarrow +\infty$. So, we compute $\sqrt{a_{ij}}$ for $a_{ij} \geq 1$ with such a precision $\epsilon$ that for any choice of $u_{ij} : |u_{ij}| \leq n$ the value of the output $\alpha = (\det B)^2$ gets computed with an error at most $10^{-n}$ (we recall that to compute the determinant we need arithmetic operations only). One can show that the bit size of $\epsilon$ can be bounded by a polynomial in the input size. We note that if $A$ is a 0-1 matrix, we don't need square roots.

The next step is to approximate sampling from the standard Gaussian distribution in $\mathbb{R}$. Let us choose a variation of the "polar method" (see Section 3.4.1 (C) of [17]). It can be shown that if $x$ and $y$ are independent random variables, uniformly distributed on $[0, 1]$, then $u = \sin(2\pi y)\sqrt{-2\ln x}$ has the standard Gaussian distribution in $\mathbb{R}$. Let $(x_{ij}, y_{ij}) : 1 \leq i, j \leq n$ be the coordinates in the $2n^2$-dimensional unit cube $[0, 1]^{n^2} \times [0, 1]^{n^2} = [0, 1]^{2n^2}$ and $u_{ij}$ be the coordinates in $\mathbb{R}^{n^2}$. This allows us to construct a map

$$\Phi : [0, 1]^{n^2} \times [0, 1]^{n^2} \longrightarrow \mathbb{R}^{n^2}, \qquad \text{where} \qquad \Phi_{ij}(x, y) = \sin(2\pi y_{ij})\sqrt{-2\ln x_{ij}},$$

such that the push-forward measure of the standard Lebesgue measure $\lambda$ on $[0, 1]^{2n^2}$ is the standard Gaussian measure $\mu$ in $\mathbb{R}^{n^2}$. The output $\alpha(u_{ij})$ is a function on $\mathbb{R}^{n^2}$. The map $\Phi$ is singular on the part of the boundary where $x_{ij} = 0$ or $x_{ij} = 1$, so we find a parallelepiped $\Pi \subset [0, 1]^{2n^2}$, such that $\lambda(\Pi) \geq 0.99$ for all sufficiently large $n$ and $1 - 1/n^3 \geq x_{ij} \geq 1/n^3$ for any point in $\Pi$. The map $\Phi$ restricted to $\Pi$ is Lipschitz, so we are able to find a rational $\delta > 0$ such that if $\|z_1 - z_2\| \leq \delta$ for $z_1, z_2 \in \Pi$ then $|\alpha(\Phi(z_1)) - \alpha(\Phi(z_2))| \leq 10^{-n}$. It can be shown that the size of $\delta$ can be bounded by a polynomial in the input size. Next, we have to approximate sampling from the Lebesgue measure on $\Pi$ by the binary sampling. We do it by sampling each coordinate $x_{ij}$ and $y_{ij}$ independently. To sample a coordinate $x$ we sample $N$ random bits $b_1, \ldots, b_N$ and let

$$x = \sum_{k=1}^{N} 2^{-k} b_k.$$

The $2^{2n^2 N}$ points that can be obtained in such a way form a lattice grid in $[0,1]^{2n^2}$. We choose $N$ so large that this grid forms a $\delta$-net in $[0,1]^{2n^2}$. One can choose an $N$ that is bounded by a polynomial in $n$ and $\ln \delta$. Now we can approximate sampling from the standard Gaussian measure in $\mathbb{R}^{n^2}$: we sample a point $z \in [0,1]^{2n^2}$ as above, accept it if $z \in \Pi$ and compute $u = \Phi(a)$ approximately with a sufficiently high precision, so that the corresponding value of $\alpha(u) = (\det B)^2$ gets computed with an error at most $2 \cdot 10^{-n}$ (taking into account the error from approximate computation of $\sqrt{a_{ij}}$). Again, this can be done in polynomial time.

Finally, the output $\alpha \leq 3 \cdot 10^{-n}$ is rounded to 0.

It is interesting to note, that even to approximate the permanent of a 0-1 matrix (a problem, which sounds purely combinatorial) we seem to have to deal with approximate computation of such "non-combinatorial" functions as sin and ln.

Algorithm 2.3 is modified similarly, except that on the first step we perturb the matrices $Q_i \longmapsto Q_i + \epsilon I$, where $I$ is the identity matrix and $\epsilon > 0$ is such that the value of $D(Q_1, \ldots, Q_n)$ changes by not more than $10^{-n}$ (the bit size of $\epsilon$ can be bounded by a polynomial in the size of the input). Now $Q_i$ are strictly positive definite and this makes computation of the decomposition $Q_i = T_i T_i^*$ stable, so that we can compute $T_i$ with the desired precision.

**(5.2) An application of the mixed discriminant to counting.** Suppose we are given a rectangular $n \times m$ matrix $A$ with the columns $u_1, \ldots, u_m$, which we interpret as vectors from $\mathbb{R}^n$. Suppose that for any subset $I \subset \{1, \ldots, m\}$, $I = \{i_1, \ldots, i_n\}$, the determinant of the submatrix $A_I = [u_{i_1}, \ldots, u_{i_n}]$ is either 0, $-1$ or 1. Such an $A$ represents a *unimodular* matroid on the set $\{1, \ldots, m\}$ and a subset $I$ with $\det A_I \neq 0$ are called a *base* of the matroid (see [22]).

Suppose now, that the columns of $A$ are colored into $n$ different colors, which induces a partition $\{1, \ldots, m\} = J_1 \cup \ldots \cup J_n$. We are interested in the number of bases that have precisely one index of each color. Let us define the positive semidefinite matrices $Q_1, \ldots, Q_n$ as follows:

$$Q_k = \sum_{i \in J_k} u_i \otimes u_i, \quad k = 1, \ldots, n.$$

Then the number of such bases can be expressed as $D(Q_1, \ldots, Q_n)$. Indeed, using the linearity of the mixed discriminant (Section 1.2) and Lemma 3.4, we have

$$D(Q_1, \ldots, Q_n) = \sum_{I = \{i_1, \ldots, i_n\}} D(u_{i_1} \otimes u_{i_1}, \ldots, u_{i_n} \otimes u_{i_n})$$

$$= \sum_{I = \{i_1, \ldots, i_n\}} \left( \det[u_{i_1}, \ldots, u_{i_n}] \right)^2,$$

where the sums are taken over all $n$-subsets $I$, having precisely one element of each color and the proof follows.

**(5.2.1) Example. Trees in a graph.** Suppose we have a connected graph $G$ with $n$ vertices and $m$ edges. Suppose further, that the edges of $G$ are colored into $n - 1$ different colors. We are interested in spanning trees $T$ in $G$ such that all edges of $T$ have different colors. Let us number the vertices of $G$ by $1, \ldots, n$ and the edges of $G$ by $1, \ldots, m$. Let us make $G$ an oriented graph by orienting its edges arbitrarily. We consider the truncated *incidence matrix* (with the last row removed) $A = (a_{ij})$ for $1 \leq i \leq n - 1$ and $1 \leq j \leq m$ as an $(n - 1) \times m$ matrix such that

$$a_{ij} = \left\{ \begin{array}{rl} 1 & \text{if } i \text{ is the beginning of } j \\ -1 & \text{if } i \text{ is the end of } j \\ 0 & \text{otherwise.} \end{array} \right.$$

The spanning trees of $G$ are in one-to-one correspondence with non-degenerate $(n - 1) \times (n - 1)$ submatrices of $A$ and the determinant of a such a submatrix is either $1$ or $-1$ (see, for example, Chapter 4 of [7]). Hence counting colored trees reduces to computation of the mixed discriminant of some positive semidefinite matrices, computed from the incidence matrix of the graph.

## Acknowledgment

This research was partially supported by Alfred P. Sloan Research Fellowship, by NSF grant DMS 9501129 and by the Mathematical Sciences Research Institute, Berkeley, CA through NSF grant DMS 9022140.

Alexander Barvinok, Department of Mathematics, University of Michigan, Ann Ar-
bor, MI 48109-1109, USA
  *E-mail address*: barvinok@math.lsa.umich.edu